\documentstyle[12pt, amsfonts]{amsart}
\def\H{{\mathbb H}}
\def\Cnh{{\C^{n,1}}}
\def\proj{{\mathbb{P}}}
\def\chs{{\H_\C^n}}
\def\C{{\mathbb{C}}}

\def\Rn{{\R^n}}
\def\N{{\mathbb{N}}}
\def\s{{\mathrm{S}}}
\def\sn{{\s^{n-1}}}
\def\B{{\mathrm{B}}}
\def\Vol{{\mathrm{Vol}}}
\def\HVol{{\mathrm{HVol}}}

\def\sw{{\mathcal S}}

\def\note#1{\ifvmode\leavevmode\fi\vadjust{\vbox to0pt{\vss
 \hbox to 0pt{\hskip\hsize\hskip1em
\vbox{\hsize2cm\small\raggedright\pretolerance10000
 \noindent #1\hfill}\hss}\vbox to8pt{\vfil}\vss}}}
 
\renewcommand{\:}{\, : \,}

\begin{document}
\def\R{{\mathbb R}}
\def\Z{{\mathbb Z}}
\def\C{{\mathbb C}}
\newcommand{\trace}{\rm trace}
\newcommand{\Ex}{{\mathbb{E}}}
\newcommand{\Prob}{{\mathbb{P}}}
\newcommand{\E}{{\cal E}}
\newcommand{\F}{{\cal F}}
\newtheorem{df}{Definition}
\newtheorem{theorem}{Theorem}
\newtheorem{lemma}{Lemma}
\newtheorem{pr}{Proposition}
\newtheorem{co}{Corollary}
\def\sign{\mbox{ sign }}
\def\a{\alpha}
\def\N{{\mathbb N}}
\def\A{{\cal A}}
\def\L{{\cal L}}
\def\X{{\cal X}}
\def\F{{\cal F}}
\def\c{\bar{c}}
\def\diam{\mbox{\rm dim}}
\def\vol{\mbox{\rm Vol}}  
\def\b{\beta}
\def\t{\theta}
\def\l{\lambda}
\def\e{\varepsilon}
\def\colon{{:}\;}
\def\pf{\noindent {\bf Proof :  \  }}
\def\endpf{ \begin{flushright}
$ \Box $ \\
\end{flushright}}

\title[The Lower Dimensional Busemann-Petty Problem in $\chs$]{The Lower Dimensional Busemann-Petty Problem in the Complex Hyperbolic Space}
\author{Susanna Dann}
\address{Mathematics Department\\
University of Missouri\\
Columbia, MO 65211}
\email{danns@@missouri.edu}
\subjclass[2010]{42, 52}
\begin{abstract}
The lower dimensional Busemann-Petty problem asks whether origin-symmetric convex bodies in $\R^n$ with smaller volume of all $k$-di\-men\-sio\-nal sections necessarily have smaller volume. The answer is negative for $k > 3$. The problem is still open for $k=2, 3$. We study this problem in the complex hyperbolic $n$-space $\chs$ and prove that the answer is affirmative only for sections of complex dimension one and negative for sections of higher dimensions.
\end{abstract}
\maketitle

\section*{Introduction}
The Busemann-Petty problem asks the following question. Given two origin-symmetric convex bodies $K$ and $L$ in $\Rn$ such that 
$$ \Vol_{n-1}(K\cap H) \leq \Vol_{n-1}(L\cap H) $$
for every hyperplane $H$ in $\Rn$ containing the origin, does it follow that 
$$ \Vol_n(K) \leq \Vol_n(L)? $$
The answer is affirmative for $n\leq 4$ and negative for $n\geq 5$. The problem, posed in 1956 in \cite{BusemannPetty1956}, was solved in the late 90's as a result of a sequence of papers \cite{Ba,Bo,Ga1,Ga2,GardnerKoldobskySchlumprecht1999,Gi,K1,K2,LR,Lu,Pa,Zh1,Zh2}, see \cite{Koldobsky2005}, p. 3, for the history of the solution. 

Since then the Busemann-Petty problem was studied in other spaces as were its numerous generalizations. We will mention just a few examples. V. Yaskin studied the Busemann-Petty problem in real hyperbolic and spherical spaces, \cite{yaskin2006}. He showed that in the spherical space the answer is the same as in $\Rn$, but not so in the real hyperbolic space, in which case the answer is affirmative for $n\leq 2$ and negative for $n\geq 3$. A. Koldobsky, H. K{\"o}nig and M. Zymonopoulou proved in \cite{KoldobskyKonigZymonopoulou2008} that the answer to the complex version of the Busemann-Petty problem is affirmative for the complex dimension $n\leq 3$ and negative for $n\geq 4$. In \cite{dann2011} the author showed that the answer to the Busemann-Petty problem in the complex hyperbolic space is affirmative for $n\leq 2$ and negative for $n\geq 3$.

It is natural to ask what happens if hyperplane sections are replaced by sections of lower dimensions. Fix an integer $k$, $1\leq k \leq n-2$. Suppose that for every $k$-dimensional subspace $H \subset \Rn$
\begin{equation}\label{eq_lvol}
 \Vol_{k}(K\cap H) \leq \Vol_{k}(L\cap H) \, ,
\end{equation}
where $K, L$ are origin-symmetric convex bodies. Does it follow that 
$$ \Vol_n(K) \leq \Vol_n(L)? $$
This question is known as the lower dimensional Busemann-Petty problem (LDBP). For the case $k=1$ condition (\ref{eq_lvol}) means that the radius of $K$ does not exceed that of $L$ in all directions and the implication follows for all origin-symmetric star bodies. Hence for $n=3$ there is nothing to prove as $k$ can only be one. For $n=4$ the affirmative answer follows from the original Busemann-Petty problem. G. Zhang \cite{Zhang1996} proved that the answer to the LDBP with $k=n-j$ is affirmative if and only if all origin-symmetric convex bodies in $\Rn$ are generalized $j$-intersection bodies. Using this connection, J. Bourgain and G. Zhang \cite{BourgainZhang1999} established a negative answer for $3<k<n$, see also \cite{RubinZhang2004} for a corrected proof and \cite{Koldobsky2000} for a different proof of this result. The cases of two- and three-dimensional sections remain open for $n\geq 5$. The LDBP in the real hyperbolic space was studied by V. Yaskin in \cite{1yaskin2006}. Other results on the LDBP can be found in \cite{GrinbergZhang1999, milman2006, MR2349610, milman2008, MR2570666}.

In this article we consider the lower dimensional Busemann-Petty problem in the complex hyperbolic $n$-space. 
In order to define volume,  we identify $\C^n$ with $\R^{2n}$ via the mapping 
\begin{equation}\label{eq_identification}
(\xi_{11}+i\xi_{12}, \dots, \xi_{n1}+i\xi_{n2}) \mapsto (\xi_{11}, \xi_{12}, \dots, \xi_{n1}, \xi_{n2}) \, .
\end{equation}
A convex body $K$ in $\R^{2n}$ is called \textit{$R_{\theta}$-invariant}, if for every $\theta \in [0,2\pi]$ and every $\xi=(\xi_{11}, \xi_{12}, \dots, \xi_{n1}, \xi_{n2})\in \R^{2n}$
$$ \|\xi \|_K = \| R_{\theta}(\xi_{11}, \xi_{12}), \dots, R_{\theta}(\xi_{n1},  \xi_{n2}) \|_K \, ,$$
where $R_{\theta}$ stands for the counterclockwise rotation by an angle $\theta$ around the origin in $\R^2$. 

Recall that an origin-symmetric body $K$ in $\chs$ is called \textit{convex} if under the mapping (\ref{eq_identification}) it corresponds to an $R_{\theta}$-invariant body in $\R^{2n}$ contained in the open unit ball such that for any pair of points in $K\subset \R^{2n}$ the geodesic segment with respect to the Bergman metric on $\chs$ joining them also belongs to $K$, see \cite{dann2011} for more details. Bodies in $\R^{2n}$ contained in the open unit ball and satisfying the latter condition will be called \textit{h-convex}. We denote the volume element on $\chs$ by $d\mu_n$ and the volume of a body $K$ in $\R^{2n}$ with respect to this volume element by $\HVol_{2n}(K)$ to distinguish from the Euclidean volume of $K$.

Now the lower dimensional Busemann-Petty problem in $\chs$ can be posed as follows. Let $1\leq k \leq n-2$ and let $K$ and $L$ be two $R_{\theta}$-invariant $h$-convex bodies in $\R^{2n}$ such that 
$$ \HVol_{2k}(K\cap H) \leq \HVol_{2k}(L\cap H) $$
for every complex linear subspace $H$ of complex dimension $k$, does it follow that 
$$ \HVol_{2n}(K) \leq \HVol_{2n}(L)? $$
In this paper we prove that the answer to this problem is affirmative only for $k=1$ and negative for $2\leq k\leq n-2$.

\section{Preliminaries}\label{SectionPreliminaries}
\subsection{Complex Hyperbolic Space}
We will work with the ball model of the complex hyperbolic space. The material of sections \ref{subsubsectionTheBallModel} and \ref{subsubsectionTheBergmanMetric} is taken from the book by Goldman \cite{goldman99}. We refer the interested reader to this book for more information. 
\subsubsection{The Ball Model}\label{subsubsectionTheBallModel}
Let $V$ be a complex vector space. The \textit{projective space associated to $V$} is the space $\proj(V)$ of all lines in $V$, i.e. one dimensional complex linear subspaces through the origin. 

Let $\C^{n,1}$ be the $(n+1)$-dimensional complex vector space consisting of $(n+1)$-tuples 
$$ Z = \left[ \begin{array}{l} 
									Z' \\
             			Z_{n+1}         
             	\end{array} \right] \in \C^{n+1}	$$ 
and equipped with the indefinite \footnote{neither positive- nor negative-semidefinite} Hermitian form 
\begin{align*}
     \left\langle Z, W \right\rangle 	&:= (Z', W')- Z_{n+1} \overline{W}_{n+1} \\
     																	&= Z_1 \overline{W}_1 + \dots + Z_n \overline{W}_n - Z_{n+1} \overline{W}_{n+1} \, ,
\end{align*}
where $Z'$ is a vector in $\C^n$ and $Z_{n+1}\in\C$. Consider the subset of \textit{negative vectors in $\C^{n,1}$}
$$ N := \{Z\in\C^{n,1} \:  \left\langle Z,Z \right\rangle <0 \} .$$
The \textit{complex hyperbolic n-space $\chs$} is defined to be $\proj(N)$, i.e. the subset of $\proj( \Cnh )$ consisting of negative lines in $\C^{n,1}$. We identify $\chs$ with the open unit ball 
$$ \B^n := \{ z\in\C^n \: (z,z) < 1 \} $$
as follows. Define a mapping $A$ by 
\begin{align*}
     A &: \C^n \longrightarrow \proj(\C^{n,1}) \\
       & z'\longmapsto \left[ \begin{array}{l} z' \\	1	\end{array} \right] \,.
\end{align*}
Since for negative vectors in $\C^{n,1}$ the $(n+1)$-coordinate is necessarily different from zero, $\chs \subset A(\C^n)$. The mapping $A$ identifies $\B^n$ with $\chs$ and $\partial \B^n = \s^{2n-1} \subset \C^n $ with $\partial \chs$.
\begin{theorem}\textnormal{(\cite{goldman99}, Theorem 3.1.10)}
Let $F\subset \proj (\C^{n,1})$ be a complex $m$-di\-men\-sional projective subspace which intersects $\chs$. Then $F \cap \chs$ is a totally geodesic holomorphic submanifold biholomorphically\footnote{biholomorphic mapping = conformal mapping} isometric to $\H_\C^m$.
\end{theorem}
\noindent
The intersection of $\chs$ with a complex hyperplane is a totally geodesic holomorphic complex hypersurface, called a \textit{complex hyperplane in $\chs$}. Its boundary is a smoothly embedded $(2n-3)$-sphere in $\partial \chs$. 

\subsubsection{The Bergman Metric and the Volume Element}\label{subsubsectionTheBergmanMetric}
We normalize the \textit{Berg\-man metric}, a Hermitian metric on $\chs$, to have constant holomorphic sectional curvature $-1$. It can be described as follows. Let $x,y$ be a pair of distinct points in $\B^n$ and let $\overleftrightarrow{xy}$ denote the unique complex line they span. The Bergman metric restricts on $\overleftrightarrow{xy}\cap\B^n$ to the Poincar{\'e} metric of constant curvature $-1$ given by: 
$$ \frac{4 R^2 dz d\overline{z}}{(R^2-r^2)^2} \, ,$$
where $R$ is the radius of the disc $\overleftrightarrow{xy}\cap\B^n$ and $r=r(z)$ is the Euclidean distance to the center of the disc $\overleftrightarrow{xy}\cap\B^n$. As $\overleftrightarrow{xy}$ is totally geodesic, the distance between $x$ and $y$ in $\chs$ equals the distance between $x$ and $y$ in $\overleftrightarrow{xy}\cap\B^n$ with respect to the above Poincar{\'e} metric. Moreover, the geodesic from $x$ to $y$ in $\chs$ is the Poincar{\'e} geodesic in $\overleftrightarrow{xy}\cap\B^n$ joining $x$ and $y$. The Poincar{\'e} geodesics are circular arcs orthogonal to the boundary and straight lines through the center. 

The volume element on $\chs$ is 
$$d\mu_n = 8^n \frac{r^{2n-1} dr d\sigma}{(1-r^2)^{n+1}} $$
where $d\sigma$ is the volume element on the unit sphere $\s^{2n-1}$. Thus for a subset $K$ of $\chs$ we have 
\begin{equation*}
\HVol_{2n}(K) = \int_{K} d\mu_{n} = 8^{n} \int_{\s^{2n-1}} \int_0^{\|x\|_K^{-1}} \frac{r^{2n-1} dr d\sigma}{(1-r^2)^{n+1}} \, .
\end{equation*} 
Let $H\subset \C^n$ be a complex subspace of dimension $n-k$: $\dim_{\C}(H)=n-k$, then
\begin{equation*}
\HVol_{2n-2k}(K\cap H) = \int_{K\cap H} d\mu_{n-k} = 8^{n-k} \int_{\s^{2n-1}\cap H} \int_0^{\|x\|_K^{-1}} \frac{r^{2n-2k-1} dr d\sigma}{(1-r^2)^{n-k+1}} \, .
\end{equation*}  

\subsubsection{Origin Symmetric Convex Sets in $\chs$}\label{subsubsection_oscs}
Origin symmetric convex bodies in $\chs$ are $R_{\theta}$-invariant bodies in $\R^{2n}$ contained in the open unit ball that are geodesically convex with respect to the Bergman metric. While it is not true in general that a convex body contained in the open unit ball is $h$-convex, see \cite{dann2011} for a counterexample, one can dilate a convex body of strictly positive curvature to make it $h$-convex. 

\begin{lemma}\label{lemma_hconvex}\textnormal{(\cite{dann2011}, Lemma 1)}
Let $D$ be an origin-symmetric convex body in $\R^{2n}$ of strictly positive curvature. Then there is an $\alpha>0$ so that the dilated body $\alpha D$ is $h$-convex.
\end{lemma}

Recall that any real ellipsoid in $\R^n$ has a section by a two-dimensional plane that is a circle. Moreover, two-dimensional planes parallel to the circular section also intersect the ellipsoid in a circle. These facts seem to belong to the folklore of the theory of the second-order surfaces, see, for instance, p. 17-18 in \cite{MR0046650} (we give proofs of these facts in the appendix for the convenience of the reader). We call a \text{\it complex ellipsoid} a real $R_{\theta}$-invariant ellipsoid in $\R^{2n}$. 

\begin{theorem}\label{th_ce}
Complex ellipsoids contained in the open unit ball are $h$-convex.
\end{theorem}

\pf
The $R_{\theta}$-invariance implies that all sections of a complex ellipsoid by one-dimensional complex subspaces are circles. Hence, all non-empty sections by affine one-dimensional complex planes are circles as well. Such circles are geodesically convex, since the geodesics are circular arcs orthogonal to the boundary of the ball.  
\endpf

Let us make an observation that will be used later. Consider a disc obtained as the intersection of the open unit ball with an arbitrary affine complex plane and a circle of radius $r<1$ centered at the origin. Move this circle, inside the disc, in any fixed direction. Ones the circle starts intersecting the boundary of the disc, consider only the arc of this circle that lies inside the disc. The question is: How far along this fixed direction can one move the circle and preserve its geodesic convexity? The geodesic convexity of the circle will be preserved all the way till the circle intersects the boundary of the ball orthogonally.

\subsection{Convex Geometry}

\subsubsection{Basic Definitions}
The main tool used in this paper is the Fourier transform of distributions, see \cite{GelfandShilov1964} as the classical reference for this topic. As usual, denote by $\sw(\Rn)$ the \textit{Schwartz space} of rapidly decreasing infinitely differentiable functions on $\Rn$, also referred to as \textit{test functions}, and by $\sw'(\Rn)$ the space of \textit{distributions} on $\Rn$, the continuous dual of $\sw(\Rn)$. The Fourier transform $\hat{f}$ of a distribution $f$ is defined by $\left\langle \hat{f}, \varphi \right\rangle = \left\langle f, \hat{\varphi} \right\rangle$ for every test function $\varphi$. A distribution $f$ on $\Rn$ is \textit{even homogeneous of degree $p\in\R$}, if 
$$ \left\langle f(x), \varphi\left(\frac{x}{\alpha}\right) \right\rangle = |\alpha|^{n+p} \left\langle f, \varphi  \right\rangle $$
for every test function $\varphi$ and every $\alpha\in\R, \alpha\neq 0$. The Fourier transform of an even homogeneous distribution of degree $p$ is an even homogeneous distribution of degree $-n-p$. We call a distribution $f$ \textit{positive definite}, if for every test function $\varphi$
$$ \left\langle f(x), \varphi \ast \overline{\varphi}(-x) \right\rangle \geq 0 \, . $$  
This is equivalent to $\hat{f}$ being a positive distribution, i.e. $\left\langle \hat{f}, \varphi \right\rangle \geq 0$ for every non-negative test function $\varphi$.

A compact subset $K$ of $\R^n$ containing the origin as an interior point is called a \textit{star body} if every line through the origin crosses the boundary in exactly two points different from the origin, and its \textit{Minkowski functional} is defined by
$$ \|x\|_K := \min \{ a\geq 0 \: x\in aK\} \, .$$  
The boundary of $K$ is continuous if $\|\cdot \|_K$ is a continuous function on $\Rn$. If in addition $K$ is origin-symmetric and convex, then the Minkowski functional is a norm on $\R^n$. A star body $K$ is said to be \textit{$k$-smooth}, $k\in\N\cup \{0\}$, if the restriction of $\|\cdot \|_K$ to the unit sphere $\sn$ belongs to the class $C^k(\sn)$ of $k$-times continuously differentiable functions on the unit sphere. If $\|\cdot \|_K \in C^k(\sn)$ for any $k\in\N$, then a star body $K$ is said to be \textit{infinitely smooth}. For $x\in\sn$, the \textit{radial function of $K$},  $\rho_K(x)=\|x\|_K^{-1}$, is the Euclidean distance from the origin to the boundary of $K$ in the direction $x$. All bodies considered in the sequel contain the origin as an interior point. 

\subsubsection{Fourier Approach to Sections}\label{sssec_FAtS}
It was shown in \cite{Koldobsky2005}, Lemma 3.16, that for an infinitely smooth origin-sym\-me\-tric star body $K$ in $\Rn$ and $0<p<n$, the Fourier transform of the distribution $\|x\|^{-p}_K$ is an infinitely smooth function on $\Rn \setminus \{0\}$, homogeneous of degree $-n+p$. We shall use a version of the Parseval's formula on the sphere:

\begin{lemma}\textnormal{(\cite{Koldobsky2005}, Lemma 3.22)}\label{lemma_ParsevalOnTheSphere}
Let $K$ and $L$ be infinitely smooth origin-symmetric star bodies in $\R^n$, and let $0<p<n$. Then 
$$ \int_{\sn} (\|\cdot\|_K^{-p})^{\wedge}(\theta) (\|\cdot\|_L^{-n+p})^{\wedge}(\theta) d\theta = (2\pi)^n \int_{\sn} \|\theta\|_K^{-p} \|\theta\|_L^{-n+p} d\theta \, . $$
\end{lemma} 

The classes of $k$-intersection bodies were introduced by A. Koldobsky in \cite{Koldobsky1999, Koldobsky2000} as follows. Let $1\leq k < n$ and let $D$ and $L$ be origin-symmetric star bodies in $\R^n$. The body $D$ is called a \textit{$k$-intersection body of $L$} if for every $(n-k)$-dimensional subspace $H$ of $\R^n$
$$ \Vol_k(D \cap H^{\perp}) = \Vol_{n-k}(L \cap H) \, .$$
An origin-symmetric star body $K$ in $\Rn$ is a $k$-intersection body if and only if $\|\cdot \|_K^{-k}$ is a positive definite distribution on $\Rn$. 

Let $0<k<n$ and let $H$ be an $(n-k)$-dimensional subspace of $\Rn$. Fix an orthonormal basis $e_1, \dots, e_k$ in the orthogonal subspace $H^{\perp}$. For a star body $K$ in $\Rn$, define the \textit{$(n-k)$-dimensional parallel section function $A_{K,H}$} as a function on $\R^k$ such that for $u\in \R^k$
\begin{align*}
	A_{K,H}(u)	&= \Vol_{n-k}(K\cap\{H+u_1 e_1 + \cdots + u_k e_k\}) \\
       				&= \int_{\{ x\in \Rn\: (x,e_1)=u_1, \dots , (x,e_k)=u_k \}} \chi (\|x\|_K) dx\, ,
\end{align*}
where $\chi$ is the indicator function of the interval $[0,1]$. If $K$ is infinitely smooth, the function $A_{K,H}$ is infinitely differentiable at the origin. We shall make use of the following fact:
\begin{lemma}\label{lem_PSFE}\textnormal{(\cite{Koldobsky2000}, Theorem 2)}
Let $K$ be an infinitely smooth origin-symmetric star body in $\Rn$ and $0<k<n$. Then for every $(n-k)$-dimensional subspace $H$ of $\Rn$ and for every $m\in \N\cup\{0\}$, $m<(n-k)/2$,
$$ \Delta^m A_{K,H}(0) = \frac{(-1)^m}{(2\pi)^k(n-2m-k)} \int_{\sn\cap H^{\perp}} (\|x\|_K^{-n+2m+k})^{\wedge}(\xi) d\xi \, ,$$
where $\Delta$ denotes the Laplacian on $\R^k$.
\end{lemma}
\subsubsection{Approximation Results}\label{subsubsection_AR}
One can approximate any convex body $K$ in $\Rn$ from inside or from outside in the \textit{radial metric} 
$$ \rho(K,L):= \max\limits_{x\in\sn} |\rho_K(x)-\rho_L(x)| $$
by a sequence of infinitely smooth convex bodies with the same symmetries as $K$, see Theorem 3.3.1 in \cite{schneider1993}. In particular, any $R_{\theta}$-invariant convex body in $\R^{2n}$ can be approximated by infinitely smooth $R_{\theta}$-invariant convex bodies. Any $k$-smooth star body $K$ can be approximated by a sequence of infinitely smooth star bodies $K_m$ so that the radial functions $\rho_{K_m}$ converge to $\rho_K$ in the metric of the space $C^k(\sn)$, see \cite{Koldobsky2005}, p. 27, preserving the symmetries of $K$ as well.

A convex body can also be approximated in the radial metric by convex bodies of strictly positive curvature. We shall use the following lemma from \cite{Koldobsky2005}:
\begin{lemma}\textnormal{(\cite{Koldobsky2005}, Lemma 4.10)}\label{lemma_Appr}
Let $1\leq k < n$. Suppose that $D$ is an origin-symmetric convex body in $\R^n$ that is not a $k$-intersection body. Then there exists a sequence $D_m$ of origin-symmetric convex bodies so that $D_m$ converges to $D$ in the radial metric, each $D_m$ is infinitely smooth, has strictly positive curvature and each $D_m$ is not a $k$-intersection body.
\end{lemma}

\noindent
Moreover, if $D$ is $R_{\theta}$-invariant, one can choose $D_m$ with the same property. The proof of the above lemma is based on the following fact, which allows for more general approximation results.  

\begin{lemma}\textnormal{(\cite{Koldobsky2005}, Lemma 3.11, (i))}\label{lemma_GenAppr}
Suppose that $p>-n$ and let $f_k$, $k\in\N$, and $f$ be even continuous functions on the sphere $\s^{n-1}$ so that $f_k\rightarrow f$ in $C(\s^{n-1})$. Then for every even test function $\phi$
$$ \lim\limits_{k \to \infty} \left\langle (f_k(\theta)r^p)^{\wedge}, \phi\right\rangle = \left\langle (f(\theta)r^p)^{\wedge}, \phi\right\rangle \, .$$
\end{lemma}

Next lemma translates the $R_{\theta}$-invariance of a body $K$ into a certain invariance of the Fourier transform of its Minkowski functional.

\begin{lemma}\label{lem_ConD}\textnormal{(\cite{KoldobskyKonigZymonopoulou2008}, Lemma 3)}
Suppose that $K$ is an infinitely smooth $R_{\theta}$-in\-va\-ri\-ant star body in $\R^{2n}$. Then for every $0<p<2n$ and $\xi\in\s^{2n-1}$ the Fourier transform of the distribution $\|x\|_K^{-p}$ is a constant function on $\s^{2n-1}\cap H_{\xi}^{\perp}$.
\end{lemma}

An important question is the following. For what $p$, $0<p<2n$, does the space $(\R^{2n}, \|\cdot\|_K)$ embed in $L_{-p}$, where $K$ is an $R_{\theta}$-invariant convex body? It was answered in \cite{KoldobskyKonigZymonopoulou2008}.

\begin{theorem}\label{th_IntersectionBodies}\textnormal{(\cite{KoldobskyKonigZymonopoulou2008}, Theorem 3)}
Let $n\geq 3$. Every $R_{\theta}$-invariant convex body $K$ in $\R^{2n}$ is a $(2n-4)$-, $(2n-3)$-, $(2n-2)$-, and $(2n-1)$-intersection body. Moreover, the space $(\R^{2n}, \|\cdot\|_K)$ embeds in $L_{-p}$ for every $p\in[2n-4, 2n)$.
For $n=2$, the space $(\R^{2n}, \|\cdot\|_K)$ embeds in $L_{-p}$ for every $p\in(0, 4)$.
\end{theorem}

There are examples of origin-symmetric $R_{\theta}$-invariant convex bodies in $\R^{2n}$, $n \geq 3$, that are not $k$-intersection bodies for any $1\leq k < 2n-4$. Denote by $B^n_q$ the unit ball of the complex space $l^n_q$ considered as a subset of $\R^{2n}$:
$$ B_q^n=\{\xi \in \R^{2n} \: \|\xi\|_q=((\xi_{11}^2 + \xi_{12}^2)^{q/2} + \cdots + (\xi_{n1}^2 + \xi_{n2}^2)^{q/2})^{1/q} \leq 1 \}\, . $$
For $q\geq 1$, $B_q^n$ is an origin-symmetric $R_{\theta}$-invariant convex body in $\R^{2n}$.

\begin{theorem}\label{th_ExampleOfNotIntersectionBody}\textnormal{(\cite{KoldobskyKonigZymonopoulou2008}, Theorem 4)}
For $q >2$ the space $(\R^{2n}, \|\cdot\|_q)$ does not embed in $L_{-p}$ with $0 < p < 2n-4$. In particular, the body $B_q^n$ is not a $k$-intersection body for any $1\leq k < 2n-4$. 
\end{theorem}

The following proposition is an analog of Lemma 3.4 in \cite{1yaskin2006} for $R_{\theta}$-invariant functions.

\begin{pr}\label{pr_ExistanceOfFunction}
Let $l$ be an integer, $2\leq 2l \leq 2n-2$. Let $f$ be an infinitely differentiable $R_{\theta}$-invariant function on $\s^{2n-1}$ so that its homogeneous extension of degree $-2l$ to $\R^{2n}$ is not a positive definite distribution. Then there exists an infinitely differentiable $R_{\theta}$-invariant function $g$ on $\s^{2n-1}$ such that
\begin{equation*}
\int_{\s^{2n-1}} f(x) g(x) dx >0 \, ,
\end{equation*}  
and for any $(n-l)$-dimensional complex subspace $H$ of $\C^n$
\begin{equation*}
\int_{\s^{2n-1}\cap H} g(x) dx \leq 0 \, .
\end{equation*}  
\end{pr}

\pf
Since the function $f$ is infinitely differentiable on $\s^{2n-1}$, the Fourier transform of its homogeneous extension $\left(f\left(\frac{x}{|x|} \right) |x|^{-2l} \right)^{\wedge}$ is a continuous function on $\R^{2n}\setminus \{0\}$, see Section \ref{sssec_FAtS}. By our hypothesis there exists an element on the sphere $\xi\in \s^{2n-1}$ with $\left(f\left(\frac{x}{|x|} \right) |x|^{-2l} \right)^{\wedge}(\xi) <0$ and hence, by continuity, an open subset $\Omega\subset \s^{2n-1}$ on which this function is negative. Moreover, $R_{\theta}$-invariance of the function implies $R_{\theta}$-invariance of the set $\Omega$. Choose a non-positive infinitely differentiable $R_{\theta}$-invariant function $h$ supported in $\Omega$ and extend it to a homogeneous function of degree $-2l$ on $\R^{2n}$. The Fourier transform of this extension is an $R_{\theta}$-invariant infinitely differentiable function on $\R^{2n}\setminus \{0\}$, homogeneous of degree $-2n+2l$, that is
\begin{equation*}
\left(h\left(\frac{x}{|x|} \right) |x|^{-2l} \right)^{\wedge}(y) = g\left(\frac{y}{|y|} \right) |y|^{-2n+2l}\, ,
\end{equation*}  
for some infinitely differentiable $R_{\theta}$-invariant function $g$ on $\s^{2n-1}$. $g$ is the function we seek. To see this, we compute, applying Parseval's formula on the sphere:
\begin{align*}
	\int\limits_{\s^{2n-1}} f(x) g(x) dx 
		&= \int\limits_{\s^{2n-1}} \left( f\left(\frac{x}{|x|} \right) |x|^{-2l} \right) \left(g\left(\frac{x}{|x|} \right) |x|^{-2n+2l} \right) dx \\
		&= \frac{1}{(2\pi)^{2n}}	\!\!\!\int\limits_{\s^{2n-1}} \!\!\! \left( f\left(\frac{x}{|x|} \right) |x|^{-2l} \right)^{\wedge}\!\!\!(\theta) \left(g\left(\frac{x}{|x|} \right) |x|^{-2n+2l} \right)^{\wedge}\!\!\!(\theta) d\theta \\
		&= \int\limits_{\s^{2n-1}} \left( f\left(\frac{x}{|x|} \right) |x|^{-2l} \right)^{\wedge}\!\!\!(\theta) h(\theta) d\theta >0\, ,
\end{align*}
since $h$ is non-positive with support in the set where $\left(f\left(\frac{x}{|x|} \right) |x|^{-2l} \right)^{\wedge}$ is negative.
By \cite{KoldobskyPaourisZymonopoulou2011} Proposition 4, we compute further 
\begin{align*}
	(2\pi)^{2l}\int_{\s^{2n-1}\cap H} g(x) dx 
 		&= \int_{\s^{2n-1}\cap H^{\perp}} \left(g\left(\frac{x}{|x|} \right) |x|^{-2n+2l}\right)^{\wedge}(\theta) d\theta \\
 		&= (2\pi)^{2n} \int_{\s^{2n-1}\cap H^{\perp}} h(\theta) d\theta \leq 0 \,.
\end{align*}
\endpf

\section{Solution of the Problem}

In $\H_\C^1$ all $R_{\theta}$-invariant bodies are discs and in $\H_\C^2$ one can only consider one-dimensional sections, which corresponds to the Busemann-Petty problem. Hence the LDBP in $\chs$ makes sense for $n\geq 3$. From now on we assume that $n\geq 3$.

\begin{theorem}\label{th_main}
Let $l$ be an integer with $1\leq l \leq n-2$. Then there are $R_{\theta}$-invariant $h$-convex bodies $K, L$ in $R^{2n}$ so that for every $(n-l)$-dimensional complex subspace $H$
$$ \HVol_{2n-2l}(K\cap H)\leq \HVol_{2n-2l}(L\cap H) \, , $$
but
$$ \HVol_{2n}(K) > \HVol_{2n}(L) \, . $$
\end{theorem}

\noindent
\text{\bf Proof of the case $1\leq l \leq n-3$\,:\,\,\,}
Let $M$ be an $R_{\theta}$-invariant convex body in $\R^{2n}$ for which the distribution $\|\cdot\|^{-2l}_M$ is not positive definite. For example take the unit ball $B^n_q$ with $q>2$ of the complex space $l^n_q$, see Theorem \ref{th_ExampleOfNotIntersectionBody} above. There is a sequence of infinitely smooth $R_{\theta}$-invariant convex bodies $M_j$ of strictly positive curvature converging to the body $M$ in the radial metric so that the corresponding distributions $\|\cdot\|^{-2l}_{M_j}$ are not positive definite, see Lemma \ref{lemma_Appr}. Pick a body $M_j$ in this approximating sequence. 

We can dilate $M_j$ to make it $h$-convex, see Lemma \ref{lemma_hconvex}. This fact exploits the idea that locally Riemannian manifolds are close to being Euclidean. Indeed, pick a small neighborhood of the origin in $B^n$, say a ball of radius $r$. It is not hard to show, see Lemma 1 in \cite{dann2011} for details, that the Euclidean curvature of all geodesics with respect to the Bergman metric in this neighborhood is less than $\frac{2r}{1-r^2}$. On the other hand, if the smallest normal curvature among the points on the boundary of a body is positive, say at least $d>0$, then
the curvature of any boundary curve at any point of the body dilated by $\alpha$ is at least $\frac{d}{\alpha}$. Choosing a dilation factor $\alpha$ so that the dilated body is contained in the ball of radius $r$ and $\frac{d}{\alpha}>\frac{2r}{1-r^2}$ ensures that the dilated body is $h$-convex.  

Let $\alpha$ be the dilation factor that would make the body $M_j$ $h$-convex. We dilate the body $M_j$ be a smaller dilation factor of $\frac{\alpha}{2}$ and use the letter $M$ again to denote the resulting body. This ensures that small smooth deformations of the body $M$ have positive normal curvature big enough to preserve $h$-convexity.   

Define another body $L$ by
\begin{equation}\label{eq_de}
\rho_L(\theta)=\sqrt{\frac{\rho^2_M(\theta)}{1+\rho^2_M(\theta)}} 
\end{equation}
for $\theta\in\s^{2n-1}$. Observe that the body $L$ so defined is contained in the body $M$. $L$ is infinitely smooth. The $R_{\theta}$-invariance of the body $M$ is preserved under this transformation and implies the $R_{\theta}$-invariance of $L$. The body $L$ is the image of $M$ under the smooth transformation given in polar coordinates by 
\begin{equation}\label{eq_MinPC}
(r,\theta) \mapsto \left( \sqrt{\frac{r^2}{1+r^2}}, \theta \right).
\end{equation}
This transformation preserves Euclidean convexity as it maps straight line segments to elliptic arcs. Indeed, consider a line in general position given by $y=mx+c$. Its equation in polar coordinates is $r=\frac{c}{\sin\theta-m\cos\theta}$ and its image under the map (\ref{eq_MinPC}) is $(m^2+c^2)x^2-2m xy +(1+c^2)y^2-c^2=0$, which is a general ellipse since $4m^2-4(m^2+c^2)(1+c^2)<0$.
Note that by restricting to small values of $r$, one can make $\sqrt{\frac{r^2}{1+r^2}}$ arbitrary close to $r$. Hence we can assume that the body $L$ is $h$-convex. Moreover, solving equation (\ref{eq_de}) for $\|\cdot\|_M^{-2}$, raising to power $l$ and extending by homogeneity to $\R^{2n}$, we obtain that the distribution 
$$\frac{\|x\|_L^{{-2l}}}{\left( 1-\left(\frac{|x|}{\|x\|_L}\right)^2\right)^l} = \|x\|_M^{-2l}$$
is not positive definite. 

By Proposition \ref{pr_ExistanceOfFunction} there exists an infinitely differentiable $R_{\theta}$-invariant function $g$ on $\s^{2n-1}$ satisfying 
\begin{equation}\label{eq_IntegralInequality}
\int_{\s^{2n-1}} \frac{\|x\|_L^{-2l}}{\left( 1-\|x\|^{-2}_L\right)^l} \, g(x) dx >0 \, ,
\end{equation} 
and  
\begin{equation}\label{eq_IntegralInequalityForG}
\int_{\s^{2n-1}\cap H} g(x) dx \leq 0 \, ,
\end{equation}  
for any $(n-l)$-dimensional complex subspace $H$ of $\C^n$. Define another infinitely smooth $R_{\theta}$-invariant body $K$ by 
\begin{equation}\label{eq_DefiningEquationForK}
	\int_0^{\|\theta\|^{-1}_K} \frac{r^{2n-2l-1}}{(1-r^2)^{n-l+1}} dr = \int_0^{\|\theta\|^{-1}_L} \frac{r^{2n-2l-1}}{(1-r^2)^{n-l+1}} dr + \epsilon g(\theta) \, ,
\end{equation}  
where $\theta\in \s^{2n-1}$ and $\epsilon >0$ small. Since the body $L$ has strictly positive curvature big enough to ensure its $h$-convexity, for small enough $\epsilon$ the body $K$ is also $h$-convex. This follows from essentially the same argument as for the strict convexity of small perturbations of the above form of strictly convex bodies, see for example \cite{Zvavitch2005} Proposition 2. 

Let $H$ be an $(n-l)$-dimensional complex subspace and integrate (\ref{eq_DefiningEquationForK}) over $\s^{2n-1}\cap H$,
\begin{align*}
	\int\limits_{\s^{2n-1}\cap H} \int\limits_0^{\|\theta\|^{-1}_K} \frac{r^{2n-2l-1}}{(1-r^2)^{n-l+1}} dr d\theta
	= \int\limits_{\s^{2n-1}\cap H} \int\limits_0^{\|\theta\|^{-1}_L} & \frac{r^{2n-2l-1}}{(1-r^2)^{n-l+1}} dr d\theta \\
	&+ \epsilon  \! \int\limits_{\s^{2n-1}\cap H} g(\theta) d\theta \, .
\end{align*}
By (\ref{eq_IntegralInequalityForG}) the second addend is non-positive, hence
\begin{equation*}
	\int\limits_{\s^{2n-1}\cap H} \int\limits_0^{\|\theta\|^{-1}_K} \frac{r^{2n-2l-1}}{(1-r^2)^{n-l+1}} dr d\theta \leq 	\int\limits_{\s^{2n-1}\cap H} \int\limits_0^{\|\theta\|^{-1}_L} \frac{r^{2n-2l-1}}{(1-r^2)^{n-l+1}} dr d\theta  \, ,
\end{equation*}
which means that for any $(n-l)$-dimensional complex subspace $H$
\begin{equation*}
	\HVol_{2n-2l} (K\cap H) \leq 	\HVol_{2n-2l} (L\cap H) \, .
\end{equation*} 

\noindent
Now we multiply (\ref{eq_DefiningEquationForK}) by the distribution $\frac{\|x\|_L^{{-2l}}}{\left( 1-\left(\frac{|x|}{\|x\|_L}\right)^2\right)^l}$ and integrate over $\s^{2n-1}$
\begin{align*}
	&\int\limits_{\s^{2n-1}} \frac{\|x\|_L^{-2l}}{\left( 1-\|x\|^{-2}_L\right)^l} \int\limits_0^{\|x\|^{-1}_K} \frac{r^{2n-2l-1}}{(1-r^2)^{n-l+1}} dr dx \\
	&=\int\limits_{\s^{2n-1}} \frac{\|x\|_L^{{-2l}}}{\left( 1-\|x\|^{-2}_L\right)^l} \int\limits_0^{\|x\|^{-1}_L} \frac{r^{2n-2l-1}}{(1-r^2)^{n-l+1}} dr dx + \epsilon \int\limits_{\s^{2n-1}} g(x) \frac{\|x\|_L^{{-2l}}}{\left( 1-\|x\|^{-2}_L\right)^l} dx \, .
\end{align*}  
By (\ref{eq_IntegralInequality}) the second addend in the above equality is strictly positive, hence
\begin{equation}\label{eq_InequalityForDifference}
		\int_{\s^{2n-1}} \frac{\|x\|_L^{{-2l}}}{\left( 1-\|x\|^{-2}_L\right)^l} \int_{\|x\|^{-1}_L}^{\|x\|^{-1}_K} \frac{r^{2n-2l-1}}{(1-r^2)^{n-l+1}} dr dx >0 \,.
\end{equation} 
Observe that the function $\frac{r^{2l}}{(1-r^2)^l}$ is an increasing function on the interval $(0,1)$. For $a,b\in (0,1)$
\begin{align*}
\frac{a^{2l}}{(1-a^2)^l} \int_a^b \frac{r^{2n-2l-1}}{(1-r^2)^{n-l+1}} dr 
	&= \frac{a^{2l}}{(1-a^2)^l} \int_a^b \frac{r^{2n-1}}{(1-r^2)^{n+1}} \frac{r^{-2l}}{(1-r^2)^{-l}} dr \\ 
	&= \int_a^b \frac{r^{2n-1}}{(1-r^2)^{n+1}} \frac{a^{2l}}{(1-a^2)^l} \left(\frac{r^{2l}}{(1-r^2)^{l}}\right)^{-1} dr \\
	&\leq \int_a^b \frac{r^{2n-1}}{(1-r^2)^{n+1}} dr \, ,
\end{align*} 
that is
\begin{equation}
	\frac{a^{2l}}{(1-a^2)^l} \int_a^b \frac{r^{2n-2l-1}}{(1-r^2)^{n-l+1}} dr \leq \int_a^b \frac{r^{2n-1}}{(1-r^2)^{n+1}} dr \, .
\end{equation} 
Note that this inequality is true for both $a\leq b$ and $b\leq a$. Integrating the above inequality over $\s^{2n-1}$ with $a=\|x\|^{-1}_L$ and $b=\|x\|^{-1}_K$, we get
\begin{equation*}
		\int\limits_{\s^{2n-1}} \frac{\|x\|_L^{{-2l}}}{\left( 1-\|x\|^{-2}_L \right)^l}  \int\limits_{\|x\|^{-1}_L}^{\|x\|^{-1}_K} \frac{r^{2n-2l-1}}{(1-r^2)^{n-l+1}} dr dx 
		\leq \int\limits_{\s^{2n-1}} \int\limits_{\|x\|^{-1}_L}^{\|x\|^{-1}_K} \frac{r^{2n-1}}{(1-r^2)^{n+1}} dr dx\, . 
\end{equation*} 
As the left hand side in the above inequality is strictly positive by (\ref{eq_InequalityForDifference}), the right hand side is strictly positive as well, and hence
\begin{equation*}
\int_{\s^{2n-1}} \int_0^{\|x\|^{-1}_K} \frac{r^{2n-1}}{(1-r^2)^{n+1}} dr dx > \int_{\s^{2n-1}} \int^{\|x\|^{-1}_L}_0 \frac{r^{2n-1}}{(1-r^2)^{n+1}} dr dx\, ,
\end{equation*} 
which is equivalent to
\begin{equation*}
\HVol_{2n}(K) > \HVol_{2n}(L).
\end{equation*}  
This completes the proof of the case $1\leq l \leq n-3$. 
\endpf

\noindent
\text{\bf Proof of the case $l=n-2$\,:\,\,\,}
For an element $\xi=(\xi_{11}, \xi_{12}, \cdots, \xi_{n1}, \xi_{n2})$ of $\R^{2n}$ denote by $\xi_n=(\xi_{n1}, \xi_{n2})$ and by $\tilde{\xi}=(\xi_{11}, \xi_{12}, \cdots, \xi_{(n-1)1}, \xi_{(n-1)2})$, then $\xi=(\tilde{\xi}, \xi_n)$. We will work with the following map, written in polar coordinates as 
\begin{equation}\label{eq_cem}
(r,\theta) \mapsto \left( \sqrt{\frac{r^2}{1-r^2}}, \theta \right) .
\end{equation}
This is the inverse of the map (\ref{eq_MinPC}). Note that this map, restricted to the two-dimensional plane $xy$, takes the ellipse $\frac{x^2}{a^2}+\frac{y^2}{b^2}=1$ to the curve
$$ \left(\frac{1}{a^2}-1\right)x^2+\left(\frac{1}{b^2}-1\right)y^2=1 \, ,$$
which is an ellipse for $a,b<1$, a horizontal or a vertical line for $a=1$ or $b=1$ and a hyperbola for $a>1$ or $b>1$. Indeed, writing the equation of the ellipse in polar coordinates, we obtain $\frac{r^2}{1-r^2} = \left(\left( \frac{\cos^2\theta}{a^2} + \frac{\sin^2\theta}{b^2}\right)-1 \right)^{-1}$, so the image is the curve with the equation $r^2\left( \frac{\cos^2\theta}{a^2} + \frac{\sin^2\theta}{b^2}\right)-r^2=1$. 

Let $0<s<1/2$. We will construct a body $L$ using two ellipses. Let $E_h$ be the ellipse with the axes $a=1, b=s$, denote its equation by $e_h$, and let $E_v$ be the ellipse with the axes $a=s, b>1$, with equation $e_v$. Choose $b$ close to one. As discussed above $E_h$ is mapped to the horizontal line $y=\frac{s}{\sqrt{1-s^2}}$ and $E_v$ is mapped to the hyperbola, denote it by $H$, with equation $h(y)=\pm \sqrt{\frac{s^2}{1-s^2}\left(1+\frac{b^2-1}{b^2}y^2\right)}$. Define a convex body $L$ in $\R^{2n}$ by
$$ L=\left\{\xi\in \R^{2n} \: |\tilde{\xi}| \leq e_v(|\xi_n|) \text{ and } |\xi_n|\leq e_h(|\tilde{\xi}|) \right\} \, .$$
The body $L$ as the intersection of two complex ellipsoids. The ellipsoids themselves are not entirely contained in the open unit ball, but the body $L$ is. Furthermore, by our choice of $s$ and $b$ and the observation after the Theorem \ref{th_ce} the body $L$ is $h$-convex. $L$ is not smooth along the curve where the two ellipsoids intersect, but we can make it infinitely smooth by changing it in an arbitrary small neighborhood of the boundary along this curve. Moreover, we will change $L$ in this neighborhood to a surface with a strictly positive, big enough curvature to preserve its $h$-convexity. This is possible since all ellipsoids have strictly positive curvature and the fact that by construction our ellipsoids intersect almost orthogonally; the angle of intersection approaches ninety degrees as $s$ approaches zero. Define a star body $M$ by
$$ \|x\|^{-2}_M = \frac{\|x\|_L^{{-2}}}{1-\left(\frac{|x|}{\|x\|_L}\right)^2} \, . $$ 
Since the body $L$ is contained in the unit ball, the body $M$ is well-defined and infinitely smooth. By construction both bodies are $R_{\theta}$-invariant. The body $M$ is an image of the body $L$ under the map (\ref{eq_cem}) and hence it can be described as 
$$ M=\left\{\xi\in \R^{2n} \: |\tilde{\xi}| \leq |h(|\xi_n|)| \text{ with } |\xi_n|\leq \frac{s}{\sqrt{1-s^2}} \right\} \, ,$$
except for a small neighborhood along the curve where the two complex surfaces meet. 

In dimension $2n$, Lemma \ref{lem_PSFE} with $H=H_{\xi}$, $\xi \in \s^{2n-1}$, in which case $k=2$, reads as follows: Let $K$ be an infinitely smooth origin-symmetric star body in $\R^{2n}$, then for $m\in \N\cup \{0\}$, $m<n-1$
$$ \Delta^m A_{K,H_{\xi}}(0) = \frac{(-1)^m}{(2\pi)^2 (2n-2m-2)} \int_{\s^{2n-1}\cap H_{\xi}^{\perp}} (\|x\|_K^{-2n+2m+2})^{\wedge}(\nu) d\nu \, .$$
Since the above integral is taken over the region $\s^{2n-1}\cap H_{\xi}^{\perp}$, by Lemma \ref{lem_ConD}, it follows that for an infinitely smooth $R_{\theta}$-invariant star body $K$
\begin{equation}\label{eq_PSF}
		\Delta^m A_{K,H_{\xi}}(0)=\frac{(-1)^m}{2\pi(2n-2m-2)} (\|x\|_K^{-2n+2m+2})^{\wedge}(\xi) \, .
\end{equation}
\noindent
Evaluating equation (\ref{eq_PSF}) for $m=1$, we obtain
\begin{equation}\label{eq_ourPSF}
 		\Delta A_{K,H_{\xi}}(0)=\frac{-1}{2\pi(2n-4)} (\|x\|_K^{-2n+4})^{\wedge}(\xi) \, .
\end{equation}

Let $x=(\tilde{x}, x_n)\in \R^{2n}$ with $x_n\neq (0,0)$. Choose $\xi\in \s^{2n-1}$ in the direction of $x_n$. Fix an orthonormal basis $\{e_1, e_2\}$ for $H_{\xi}^{\perp}$. For $u\in\R^2$, with $|u|<\frac{s}{2\sqrt{1-s^2}}$, compute 
\begin{align*}
	A_{M,H_{\xi}}(u) 	&= \Vol_{2n-2} (M\cap\{H_{\xi} + u_1 e_1 + u_2 e_2 \}) \\
										&= \int_{\{ x\in \R^{2n}\: (x,e_1)=u_1, (x,e_2)=u_2 \}} \chi (\|x\|_M) dx \\
										&= \int_{\s^{2n-3}} \int_0^{h(|u|)} r^{2n-3} dr d\theta \\
										&= |\s^{2n-3}| \, \frac{h(|u|)^{2n-2}}{2n-2} \\		
										&= \frac{2 \pi^{n-1}}{(n-2)!} \,  \frac{h(|u|)^{2n-2}}{2(n-1)} \\
										&= \frac{\pi^{n-1}}{(n-1)!} \, h(|u|)^{2n-2} \, ,
\end{align*}
where $|\s^{n-1}|$ stands for the surface area of the unit sphere $\sn$ in $\Rn$: $|\s^{n-1}|=2 \pi^{\frac{n}{2}}/\Gamma (\frac{n}{2})$. Thus we have
$$ A_{M,H_{\xi}}(u) = \frac{\pi^{n-1}}{(n-1)!} \left(\frac{s^2}{1-s^2}\left(1+\frac{b^2-1}{b^2}|u|^2\right)\right)^{n-1} $$
and consequently
\begin{align*}
\Delta A_{M,H_{\xi}}(u) =& \, \frac{4\pi^{n-1}}{(n-2)!} \left( \frac{s^2}{1-s^2} \right)^2 \left(\frac{s^2}{1-s^2}\left(1+\frac{b^2-1}{b^2}|u|^2\right)\right)^{n-3} \\
												& \times  \frac{b^2-1}{b^2} \left\{1+(n-1)\frac{b^2-1}{b^2}|u|^2 \right\} \, .
\end{align*}
Since $M$ is infinitely smooth we can use equation (\ref{eq_ourPSF}) to compute 
$$ (\|x\|_M^{-2n+4})^{\wedge}(\xi)= - 2\pi (2n-4)	\, \Delta A_{M,H_{\xi}}(0) = - \frac{16 \pi^n}{(n-3)!} \frac{b^2-1}{b^2}\left(\frac{s^2}{1-s^2}\right)^{n-1} . $$
As $b>1$ and $s<1$, this shows that $ \left( \frac{\|x\|_L^{{-2n+4}}}{\left(1- |x|^2 \|x\|^{-2}_L \right)^{n-2}} \right)^{\wedge}(\xi) = (\|x\|_M^{-2n+4})^{\wedge}(\xi) $ is negative in some direction $\xi$. 

Thus we have constructed an infinitely smooth $R_{\theta}$-invariant $h$-convex body $L$ of strictly positive curvature so that the distribution $\frac{\|x\|_L^{-2n+4}}{\left(1- |x|^2 \|x\|^{-2}_L \right)^{n-2}}$ is not positive definite. Now we can proceed as in the proof of the case $1\leq l \leq n-3$ to construct another infinitely smooth $R_{\theta}$-invariant $h$-convex body $K$ of strictly positive curvature so that for every $2$-dimensional complex subspace $H$
$$ \HVol_{4}(K\cap H)\leq \HVol_{4}(L\cap H) \, , $$
but
$$ \HVol_{2n}(K) > \HVol_{2n}(L) \, . $$
\endpf

Finally, observe that the case $l=n-1$ corresponds to sections by one-dimensional complex subspaces. Since sections by one-dimensional complex subspaces are discs, this case is equivalent to the case of one-dimensional sections in the real setting and so has an affirmative answer. 

Altogether we have shown that the lower dimensional Busemann-Petty problem in $\chs$ has a negative answer for sections of dimension $k\geq 2$. Case of one-dimensional sections is trivially true.  

\section*{Appendix}
Here we prove the auxiliary facts used in Section \ref{subsubsection_oscs} and some other related facts of independent interest. For illustration we consider the case of $\R^3$ first.

Let $E$ be an ellipsoid in $\R^3$. Without the loss of generality, we can assume that $E$ is centered; then its equation is
\begin{equation}\label{eq_e}
\frac{x^2}{a^2}+\frac{y^2}{b^2}+\frac{z^2}{c^2}=1 \, ,
\end{equation}
where $a>b>c$.

\begin{theorem}\label{th_csoe}
There is a section of an ellipsoid in $\R^3$ that is a circle.
\end{theorem}

\pf
Assume that the ellipsoid has the form (\ref{eq_e}). Every section of $E$ by a plane containing the $y$-axis is an ellipse having one of the axes equal to $b$. The section by the $yz$-plane is an ellipse with axes $b,c$, whereas the section by the $xy$-plane is an ellipse with axes $a,b$. Thus, by continuity, there is a plane (containing the $y$-axis) that intersects $E$ in a circle of radius $b$.

Because of the symmetry of an ellipsoid there is a second plane that intersects $E$ in a circle with radius $b$.
\endpf

Next, we would like to construct one such plane, call it $p_0$; it must have the form $\alpha x + \beta z = 0$ with $\alpha^2+\beta^2=1$. From the equation of the plane we obtain $x=-\frac{\beta}{\alpha} z$, and plug this into (\ref{eq_e}) to determine the quantity $\frac{\beta}{\alpha}$. This gives
$$ \left(\frac{\beta^2}{\alpha^2} \frac{1}{a^2}+\frac{1}{c^2}\right) z^2 + \frac{y^2}{b^2}=1 \, ,$$
which is a circle of radius $b$ if and only if
$$ \left(\frac{\beta^2}{\alpha^2} \frac{1}{a^2}+\frac{1}{c^2}\right) = \frac{1}{b^2} \text{ or equivalently } \frac{\beta^2}{\alpha^2} = \frac{a^2(c^2-b^2)}{b^2 c^2} \, .$$
Using the relation $\alpha^2+\beta^2=1$, one computes $\frac{\beta^2}{\alpha^2} = \frac{1-\alpha^2}{\alpha^2} $; this gives that
$$ \alpha^2=\frac{b^2c^2}{a^2(c^2-b^2)+b^2c^2} \text{ and }  \beta^2 = \frac{a^2(c^2-b^2)}{a^2(c^2-b^2)+b^2c^2} \, .$$ 

\begin{theorem}
Every plane parallel to the circular section of an ellipsoid in $\R^3$ also intersects the ellipsoid in a circle.
\end{theorem}

\pf
First we show it for planes parallel to $p_0$; they have the form $\alpha x+ \beta z = d$, or equivalently, $x=\frac{d-\beta z}{\alpha}$. Substituting this into (\ref{eq_e}) yields:
$$ \left(\frac{\beta^2}{\alpha^2 a^2}+\frac{1}{c^2} \right)z^2- \frac{2d \beta}{\alpha^2 a^2} z+\left(\frac{d^2}{\alpha^2 a^2}-1\right) + \frac{y^2}{b^2} =0 \, ,$$ 
writing this as $ Az^2-2Dz+Cy^2+F=0$ and completing the square, we obtain $A(z-\frac{D}{A})^2 + Cy^2 = \frac{D^2}{A}-F$, which is a circle if and only if $A=C$. But this is the same condition as for the plane $p_0$, so it is satisfied.

More generally, suppose the plane $\alpha x + \beta y + \gamma z =0$ intersects the ellipsoid in a circle. Evaluating (\ref{eq_e}) for $x=-\frac{\beta}{\alpha} y - \frac{\gamma}{\alpha} z$, yields:
\begin{equation}\label{eq_ptto} 
\left( \frac{\beta^2}{\alpha^2 a^2} + \frac{1}{b^2} \right) y^2 + 2\frac{\beta \gamma}{\alpha^2 a^2} yz + \left( \frac{\gamma^2}{\alpha^2 a^2} + \frac{1}{c^2} \right) z^2 = 1\, .
\end{equation}
To verify that parallel planes $\alpha x + \beta y + \gamma z =d$ also intersects the ellipsoid in a circle, we evaluate (\ref{eq_e}) for $x=\frac{d-\beta y -\gamma z}{\alpha}$ and obtain
$$ \left( \frac{\beta^2}{\alpha^2 a^2} + \frac{1}{b^2} \right) y^2 - 2 \frac{d \beta}{\alpha^2 a^2} y + \left( \frac{\gamma^2}{\alpha^2 a^2} + \frac{1}{c^2} \right) z^2 - 2 \frac{d \gamma}{\alpha^2 a^2} z + 2 \frac{\beta \gamma}{\alpha^2 a^2} yz = 1 - \frac{d^2}{\alpha^2 a^2}  \, .$$
Comparing this with (\ref{eq_ptto}), we see that $d$ contributes to the shift of the curve of intersection, but not to its form; it also changes the radius. Hence parallel planes that have a non-empty intersection with the ellipsoid also intersect it in a circle.
\endpf

\begin{theorem}
There is a two-dimensional section of an ellipsoid in $\R^n$ that is a circle.
\end{theorem}

\pf
We may assume that the ellipsoid has the form:
\begin{equation}\label{eq_geoae}
	\frac{x^2_1}{a^2_1} + \frac{x^2_2}{a^2_2} + \cdots + \frac{x^2_n}{a^2_n} = 1 \, ,
\end{equation}
with $a_1>a_2> \cdots > a_n$. The section with the subspace spanned by the vectors $\{x_1, x_2, x_3\}$ is an ellipsoid in $\R^3$ and the result follows by Theorem \ref{th_csoe}. 
\endpf

\begin{theorem}
Let $E$ be an ellipsoid in $\R^n$ and $H$ any subspace of $\R^n$. Sections of $E$ by affine subspaces parallel to $H$ are similar to the section by $H$ itself. 
\end{theorem}

\pf
Suppose $H$ is given by $\alpha_1 x_1 + \alpha_2 x_2 + \cdots + \alpha_n x_n = 0$ and $E$ is given by (\ref{eq_geoae}). From the equation of the plane, we have
$$ x_1=\frac{- \alpha_2 x_2 - \cdots - \alpha_n x_n}{\alpha_1} =- \frac{ \alpha_2 x_2 + y}{\alpha_1} \, .$$
Plugging this into (\ref{eq_geoae}), we obtain the equation of the section $E \cap H$:
$$ \left(\frac{\alpha_2^2}{\alpha_1^2 a_1^2} + \frac{1}{a_2^2} \right) x_2^2 + 2\frac{\alpha_2}{\alpha_1^2 a_1^2} y x_2 + \frac{1}{\alpha_1^2 a_1^2} y^2 + \frac{x^2_3}{a^2_3}+\cdots \frac{x^2_n}{a^2_n} =1 \, .$$

Affine subspaces parallel to $H$ have the form $\alpha_1 x_1 + \alpha_2 x_2 + \cdots + \alpha_n x_n = d$. From this we have 
$$ x_1=\frac{d - \alpha_2 x_2 - \cdots - \alpha_n x_n}{\alpha_1} =\frac{d - \alpha_2 x_2 - y}{\alpha_1} \, , $$
and for the section of $E$ by an affine subspace parallel to $H$, we have 
\begin{align*} 
 \left(\frac{\alpha_2^2}{\alpha_1^2 a_1^2} + \frac{1}{a_2^2} \right) x_2^2 &- 2 \frac{d \alpha_2}{\alpha^2_1 a_1^2} x_2 +2\frac{\alpha_2}{\alpha_1^2 a_1^2} y x_2 + \\
& + \frac{1}{\alpha_1^2 a_1^2} y^2 - 2 \frac{d}{\alpha^2_1 a_1^2} y
+ \frac{x^2_3}{a^2_3}+\cdots \frac{x^2_n}{a^2_n} =1-\frac{d^2}{\alpha^2_1 a_1^2} \, . 
\end{align*}

It is enough to investigate the affect of the parallel shift on any single variable, let us pick $x_2$. Comparing the above equations of the sections, we see that $d$ contributes only to the shift of the intersection surface and to the similarity factor.
\endpf

\end{document}